\documentclass[10pt]{amsart}
\usepackage{latexsym}
\usepackage{amsmath}
\usepackage{amsfonts}
\usepackage{amsthm}
\usepackage{amssymb}
\usepackage[utf8]{inputenx}
\usepackage[a5paper]{geometry}
\newcommand{\CP}[1]{\makebox{${\rm C}_p({#1})$}}
\newcommand{\USC}[1]{\makebox{${\rm USC}_p({#1})$}}
\newcommand{\Meas}[2]{\makebox{${\rm M}_p({#1},{#2})$}}
\newcommand{\USM}[2]{\makebox{${\rm USM}_p({#1},{#2})$}}
\newcommand{\seq}[2]{\left\langle #1:#2\in\omega\right\rangle}
\newcommand{\fml}[2]{\{#1:#2\in\omega\}}

\newcommand{\nula}{\makebox{${\mathbf 0}$}}
\newcommand{\jedna}{\makebox{${\mathbf 1}$}}
\newcommand{\snula}{\makebox{${\scriptstyle\mathbf 0}$}}
\newcommand{\RR}{{\mathbb R}}

\newcommand{\emm}{\bf}

\newcommand{\mc}[1]{{\mathcal #1}}
\newcommand{\Borel}{\makebox{\sc Borel}}
\newcommand{\pf}{{\parindent0cm{\bf Proof:\ }}}
\newcommand{\prfwp}{{\parindent0cm{\bf Proof\ }}}

\renewcommand{\cite}[1]{[#1]}
\newcommand{\presentingauthor}[1]{\author{#1}}
\newtheorem{theorem}{Theorem}[section]
\newtheorem{corollary}[theorem]{Corollary}
\newtheorem{lemma}[theorem]{Lemma}

\begin{document}
\title{Selection Principles for Measurable Functions  and Covering Properties}
\presentingauthor{Lev Bukovsk\'y}
\address{Institute of Mathematics, P.J. \v Saf\'arik University,}
\address{\vskip-0.6cm
Jesenn\'a 5, 041 54 Ko\v sice, Slovakia}
\thanks{This work was supported by the Grant No. 1/0097/16 of Slovak Grant Agency VEGA. The first results were presented at the conference Set Theoretic Methods in Topology and Analysis, held in the Conference Center of Institute of Mathematics of the Polish Academy of Sciences, B\c edlewo, September 3 -- 9, 2017.}
\email{lev.bukovsky{@}upjs.sk}
\date{}
\begin{abstract}
Let ${\mathcal A}\subset {\mathcal P}(X)$, $\emptyset, X\in {\mathcal A}$, ${\mathcal A}$ being closed under finite intersections. If $\psi={o},\omega,\gamma$, then
$\Psi({\mathcal A})$ is the family of those $\psi$-covers ${\mathcal U}$ for which ${\mathcal U}\subseteq {\mathcal A}$. In~\cite{BL2} I have introduced properties $(\Psi_{\snula})$ of a~family $F\subseteq {}^X\RR$ of real functions. The main result of the paper Theorem~\ref{Ufinsemi} reads as follows: if~$\Phi=\Omega,\Gamma$, then for any couple $\langle \Phi,\Psi\rangle$ different from $\langle \Omega,{\mathcal O}\rangle$, $X$ has the covering property~{\rm S}${}_1(\Phi({\mathcal A}),\Psi({\mathcal A}))$ if and only if the family of non-negative upper ${\mathcal A}$-semimeasurable real functions  satisfies the selection principle~{\rm S}${}_1(\Phi_{\snula},\Psi_{\snula})$. Similarly for {\rm S}${}_{\scriptstyle fin}$ and {\rm U}${}_{\scriptstyle fin}$. Some related results are also presented.  
 \end{abstract}
\maketitle
%%%%%%%%%%%%%%%%%%%%%%%%%%%%%%%%%%%%%%%%%%%%%%%%%%%%%%%%%%%%
%%%%%%%%%%%%%%%%%%%%%%%%%%%%%%%%%%%%%%%%%%%%%%%%%%%%%%%%%%%%
\centerline{{\rm Dedicated to Alexander Blaszczyk}}
\centerline{{ on the occasion of his 70th birthday}}
\section{Introduction}\label{Int}
One of the first results describing a~connection between covering properties of a~topological space $X$ and properties of the space of real continuous functions $\CP{X}$ was proved by W.~Hurewicz~\cite{Hu}. Indeed, in contemporary terminology, W.~Hurewicz showed that a~metric space\footnote{Actually the assertions are true for any~perfectly normal topological space, see e.g. \cite{BH}.}~$X$ is a~U${}_{\scriptstyle fin}({\mathcal O}_{\scriptstyle cnt},{\mathcal O})$-space (a~U${}_{\scriptstyle fin}({\mathcal O}_{\scriptstyle cnt},\Gamma)$-space) if and only if every continuous image of $\CP{X}$ into~${}^{\omega}\RR$ is not dominating (is eventually bounded). 
By B.~Tsaban~\cite{Ts}, a~perfectly normal topological space~$X$ has the property U${}_{\scriptstyle fin}({\mathcal O}_{\scriptstyle cnt},\Omega)$ if and only if no continuous image of $\CP{X}$ into~${}^{\omega}\RR$ is finitely dominating.
\par
In \cite{BL1} the author used upper semicontinuous functions for investigation of some covering properties.Immediately, M.~Sakai~\cite{Sa} has used this idea and has shown that a~topological space $X$ is an~S${}_1(\Gamma,\Omega)$-space if and only if the space of non-negative upper semicontinuous functions satisfies the selection property S${}_1(\Gamma_{\snula},\Omega_{\snula})$.    
\par
In \cite{BL2} we have presented a~rather systematic treatment of the connection of S${}_1$ properties of a~topological space $X$ with similar properties of sets of continuous or non-negative semicontinuous real functions defined on $X$. Now, we present general results of this type about S${}_1$, S${}_{\scriptstyle fin}$ and U${}_{\scriptstyle fin}$ that includes all the above mentioned cases. We show that similar result holds true for a~general sets of covers and corresponding families of real functions, including the Borel covers and Borel measurable functions.
%%%%%%%%%%%%%%%%%%%%%%%%%%%%%%%%%%%%%%%%%%%%%%%%%%%%%%%%%%%
%%%%%%%%%%%%%%%%%%%%%%%%%%%%%%%%%%%%%%%%%%%%%%%%%%%%%%%%%%%
\section{Terminology and Notations}\label{TN}
In the next, $X$ is an~infinite Hausdorff topological space with the topology ${\mathcal T}$ (the family of open subsets).
\par
A~family ${\mathcal U}$ of subsets of $X$ is a~{\emm cover} if $\bigcup{\mathcal U}=X$ and \mbox{$X\notin {\mathcal U}$.} From some technical reason, a~cover will be called an~$o${\emm -cover}.  A~cover~${\mathcal U}$ is {\emm large} or a~$\lambda${\emm -cover} if the set $\{U\in{\mathcal U}:x\in U\}$ is infinite for each $x\in X$. A~cover ${\mathcal U}$ is an~$\omega${\emm -cover} if for every finite set~$A\subseteq X$ there exists an~element $U\in {\mathcal U}$ such that $A\subseteq U$. An~infinite cover~${\mathcal U}$ is a~$\gamma${\emm -cover} if the set $\{U\in{\mathcal U}:x\notin U\}$ is finite for each $x\in X$. We assume that the reader knows elementary properties of introduced families of covers.
\par
We shall use the following convention. If $\varphi $ and $\psi$ denote one of the symbols $o$, $\lambda$, $\omega$ or $\gamma$, then the capital $\Phi$ and $\Psi$ denote the corresponding symbol ${\mathcal O}$, $\Lambda$, $\Omega$ or $\Gamma$, respectively. Also vice versa. 
\par
In the next we assume that the~family of sets ${\mathcal A}\subseteq {\mathcal P}(X)$ satisfies 
\begin{equation}\label{emptyX}
\emptyset, X\in {\mathcal A}. 
\end{equation}
The topology ${\mathcal T}$ on $X$, the classes ${\bold \Sigma}^0_{\xi}(X)$, ${\bold \Pi}^0_{\xi}(X)$, ${\bold \Delta}^0_{\xi}(X)$ of the Borel hierarchy, the $\sigma$-algebra $\Borel(X)$ of all Borel subsets of $X$, the $\sigma$-algebra ${\mathcal M}(X)$ of all subsets of $X$ with the Baire property or the \mbox{$\sigma$-algebra ${\mathcal L}(\RR^n)$} of all Lebesgue measurable subsets of $\RR^n$, and others, are the candidates for the family~${\mathcal A}$.
\par
If $\varphi$ is one of the symbols $o$, $\lambda$, $\omega$ or $\gamma$, then $\Phi({\mathcal A})$ is the family of all $\varphi$-covers ${\mathcal U}$ such that ${\mathcal U}\subseteq {\mathcal A}$. Finally, instead of $\Phi({\mathcal T})$ we write $\Phi(X)$  or simply $\Phi$. Evidently 
\[\Gamma({\mathcal A})\subseteq \Omega({\mathcal A})\subseteq\Lambda({\mathcal A})\subseteq {\mathcal O}({\mathcal A}).\]
\par
One can easily see, that 
\begin{eqnarray}\label{common}
{}&If\ the\ family\ {\mathcal A}\ is\ closed\ under\ finite\ intersections,\nonumber \\ 
{}&then\ for\ \Phi={\mathcal O}, \Omega,\Gamma,\ any\ two\ covers\ from\ \Phi({\mathcal A})\\ 
{}& have\ a~common\ refinement\ in\ \Phi({\mathcal A}).\nonumber
\end{eqnarray}
\par
If $\varphi$ denotes one of the symbols $o$, $\omega$ or $\gamma$, then an~open $\varphi$-cover~${\mathcal U}$ is {\emm shrinkable}, if there exists an~open~$\varphi$-cover ${\mathcal V}$ such that
\begin{equation}\label{shrink}
(\forall V\in{\mathcal V})(\exists U_V\in{\mathcal U})\,
\overline{V}\subseteq U_V.
\end{equation}
The family of all open shrinkable, open shrinkable $\omega$- and $\gamma$-covers of $X$ will be denoted by ${\mathcal O}^{\scriptstyle sh}(X)$, $\Omega^{\scriptstyle sh}(X)$ and $\Gamma^{\scriptstyle sh}(X)$, or simply ${\mathcal O}^{\scriptstyle sh}$, $\Omega^{\scriptstyle sh}$ and $\Gamma^{\scriptstyle sh}$, respectively. If $X$ is a~regular topological space, then every open cover and every open $\omega$-cover is shrinkable. 
 \par
In the next, $\varepsilon$ is always a~positive real.
\par
The set ${}^X\RR$ of all real functions on $X$ is endowed with the product topology. If $c\in \RR$ is a~real, then $\bf c$ denotes the constant function on $X$ with the value $c$. A~typical neighborhood of the function $\nula$ in the product topology on ${}^X\RR$ is the set
\begin{equation}\label{TypNZ}
N^{\varepsilon}_{x_0,\dots,x_k}=\{h\in {}^X\RR:
(\forall i\leq k)\,\vert h(x_i)\vert<\varepsilon\},
\end{equation}
where $x_0,\dots,x_k\in X$. A~sequence of~elements of ${}^X\RR$ converges to a~function $f$ in this topology, written $f_n\to f$, if it converges pointwise, i.e., if $\lim_{n\to\infty}f_n(x)=f(x)$ for each $x\in X$.
\par
If $F\subseteq {}^X\RR$ is infinite countable set, then the convergence $f_n\to f$ and the uniform convergence $f_n\rightrightarrows f$ do not depend on the enumeration of $F=\{f_n:n\in\omega\}$. So, we can say that $F$ converges to $f$ or $F$ uniformly converges to $f$, written $F\to f$ or $F\rightrightarrows f$, respectively.   
\par  
Any subset $F\subseteq {}^X\RR$ is endowed with the subspace topology, especially the subset $\CP{X}$ of all continuous real functions, simply C${}_p$, and the subset $\USC{X}$ of all upper semicontinuous real functions, simply USC${}_p$.  We denote by $F^+$ the set of all non-negative functions from a~set~$F$.\par 
If $H=\{h_0,\dots,h_k\}$ is a~finite set of real functions, then \mbox{$h=\min H$} is the function defined as
\[h(x)=\min\{h_0(x),\dots,h_k(x)\}\mbox{\ for each\ }x\in X.\]
By definition we set $\min\emptyset =\jedna$.
\par
Let ${\mathcal A}\subseteq {\mathcal P}(X)$ be a~family of subsets of $X$. A~real function $f\in{}^X\RR$ is ${\mathcal A}${\emm -measurable} if for every open interval $(a,b)$, including $a=-\infty$ and $b=+\infty$, we have $f^{-1}((a,b))\in {\mathcal A}$.  A~real function $f\in{}^X\RR$ is  {\emm  upper ${\mathcal A}$-semimeasurable} if  \mbox{$f^{-1}((-\infty,b))\in {\mathcal A}$} for every $b\in\RR$. One can easily see that if ${\mathcal A}$ is a~$\sigma$-algebra, then every~upper ${\mathcal A}$-semimeasurable function is \mbox{${\mathcal A}$-measurable.} We denote by $\Meas{X}{{\mathcal A}}$ the set of all ${\mathcal A}$-measurable real functions defined on $X$ and by $\USM{X}{{\mathcal A}}$ the set of all upper ${\mathcal A}$-semimeasurable functions defined on $X$. Both sets $\Meas{X}{{\mathcal A}}$ and $\USM{X}{{\mathcal A}}$ are endowed with the subspace topology of ${}^X\RR$. Note that 
\[\CP{X}=\Meas{X}{{\mathcal T}},\ \ \ \USC{X}=\USM{X}{{\mathcal T}}.\]
\par
If $f:X\longrightarrow {}^{\omega}\omega$, then $f$ is ${\mathcal A}$-measurable if for any set $U\subseteq {}^{\omega}\omega$, $U=\{\alpha\in{}^{\omega}\omega:(\forall k\in A)\,\alpha(k)=\beta(k)\}$, $\beta\in{}^{\omega}\omega$, $A\subseteq \omega$ finite, the set $f^{-1}(U)$ belongs to ${\mathcal A}$. 
\par
Similarly, a~function $f:X\longrightarrow {}^{\omega}\RR$ is ${\mathcal A}$-measurable if for any set \mbox{$U=\prod_{k\in\omega}U_k$,} $U_k$ is an~open interval for finitely many $k$ and $U_k=\RR$ otherwise, the set $f^{-1}(U)$ belongs to ${\mathcal A}$.
\par
We consider tree properties of a family of real functions $F\subseteq {}^X\RR$ introduced in~\cite{BL2}.
\[
\begin{array}{ll}
{}&{}\\
({\mathcal O}_{\snula})\ &\nula\notin F\textnormal{\ and\ }0\in\overline{\{f(x):f\in F\}}\textnormal{\ for every\ }x\in X.\\
{}&{}\\
(\Omega_{\snula})&\nula\notin F\textnormal{\ and\ }\nula\in\overline{F}
\textnormal{\ in the topology of\ }{}^X\RR.\\
{}&{}\\
(\Gamma_{\snula})&\nula\notin F,\ F\textnormal{\ is infinite and for every\ }x\in X\textnormal{\ and}\\
{}& \textnormal{every\ }\varepsilon\textnormal{\ the set\ }\{f\in F:\vert f(x)\vert\geq \varepsilon\}\textnormal{\ is finite.}
\\
{}&{}
\end{array}
\]
Note that a~set with the property  $(\Omega_{\snula})$ is infinite as well.
\par
Similarly as for covers we have
\[(\Gamma_{\snula})\to(\Omega_{\snula})\to({\mathcal O_{\snula}}).\]
\par
If $F\subseteq {}^X\RR$ we set
\[\vert F\vert=\{\vert f\vert:f\in F\}.\]
For $\Phi={\mathcal O},\Omega\textnormal{\ or\ }\Gamma$ one can easily see that
\begin{equation}\label{absolute}
F\ possesses\ (\Phi_{\snula})\ if\ and\ only\ if\ \vert  F\vert\ does\ so.     
\end{equation}
\par
We shall not discuss elementary properties of introduced notions. Note just that if an~infinite countable family $F$ has the property $(\Gamma_{\snula})$, then $F\to \nula$. 
\par
Let $F\subseteq {}^X\RR$ be an~infinite set of real functions. We define a~family of subsets of $X$ by setting 
\[{\mathcal U}^{\varepsilon}(F)=\{U_f^{\varepsilon}: f\in F\}
,\]
where
\begin{equation}\label{Cover}
U_f^{\varepsilon}=\{x\in X:\vert f(x)\vert<\varepsilon\}.
\end{equation}
If $f$ is an~${\mathcal A}$-measurable or non-negative upper ${\mathcal A}$-semimeasurable function then every set~$U_f^{\varepsilon}$ belongs to ${\mathcal A}$.
\par
Note that 
\begin{equation}\label{hereditary}
{\rm if\ }\varepsilon_1<\varepsilon_2,{\ then\ }
{\mathcal U}^{\varepsilon_1}(F){\ is\ a~refinement\ of\ }{\mathcal U}^{\varepsilon_2}(F).
\end{equation} 
\par
The Theorem 4.1 of \cite{BL2} reads as follows.
\begin{theorem}\label{gamma} Let $F\subseteq{}^X\RR$ be an~infinite family of real functions.  
\begin{enumerate}
\item[{\rm a)}]
The family $F$ has the property $({\mathcal O} _{\snula})$ if and only if either there exists a~subsequence of $F$ uniformly converging to~$\nula$ or there exists a~$\delta>0$ such that for every $\varepsilon<\delta$ the family ${\mathcal U}^{\varepsilon}(F)$ is an~$o$-cover of~$X$.
\item[{\rm b)}] The~family~$F$ has the property $(\Omega_{\snula})$ if and only if either there exists a~subsequence of $F$ uniformly converging to~$\nula$ or there exists \mbox{a~$\delta>0$} such that for every  $\varepsilon<\delta$ the family ${\mathcal U}^{\varepsilon}(F)$ is an~$\omega$-cover of~$X$.
\item[{\rm c)}]
The family $F$ has the property $(\Gamma_{\snula})$ if and only if either $F$ is countable and $F\rightrightarrows \nula$ or there exists $\delta>0$ such that for every $\varepsilon<\delta$ the family ${\mathcal U}^{\varepsilon}(F)$ is a~$\gamma$-cover of~$X$ and for every $f\in F$ the set $\{g\in F:U^{\varepsilon}_g=U^{\varepsilon}_f\}$ is finite.
\end{enumerate}
\end{theorem}
\par
We shall also need the Lemma~6.1 of \cite{BL2} that says:
\begin{lemma}\label{epsilontozero}
Assume that $\seq{\varepsilon_n}{n}$ is a~sequence of positive reals converging to~$0$, and \mbox{$f_n\in {}^X\RR$} for $n\in\omega$. If\/$\{U_{f_n}^{\varepsilon_n}:n\in\omega\}$ is $\lambda$-, $\omega$- or $\gamma$-cover, then either there exists an~increasing sequence $\seq{n_k}{k}$ of integers such that $f_{n_k}\rightrightarrows \nula$ or there exists a~$\delta>0$ such that for every positive $\varepsilon<\delta$, the family ${\mathcal U}^{\varepsilon}(\fml{f_n}{n})$ is an~$o$-,an~$\omega$- or a~$\gamma$-cover, respectively.
\end{lemma}
%%%%%%%%%%%%%%%%%%%%%%%%%%%%%%%%%%%%%%%%%%%%%%%%%%%%%%%%%%%
%%%%%%%%%%%%%%%%%%%%%%%%%%%%%%%%%%%%%%%%%%%%%%%%%%%%%%%%%%%
\section{Covering Propertiers}\label{Ufin}
\par
We recall some notions  introduced in~\cite{JMSS}. Let ${\mathcal A},{\mathcal B}\subseteq {\mathcal P}(Y)$ for some set $Y$.
\par
The {\emm covering property} S${}_1({\mathcal A},{\mathcal B})$ means the following: for any sequence of sets $\seq{{\mathcal U}_n\in{\mathcal A}}{n}$, for every $n$ there exists \mbox{an~$U_n\in {\mathcal U}_n$} such that $\{U_n:n\in\omega\}\in{\mathcal B}$.
\par
The {\emm covering property} S${}_{\scriptstyle fin}({\mathcal A},{\mathcal B})$ means the following: for any sequence of sets $\seq{{\mathcal U}_n\in{\mathcal A}}{n}$, for every $n$ there exists a~finite set ${\mathcal V}_n\subseteq {\mathcal U}_n$ such that $\bigcup_{n\in\omega}{\mathcal V}_n\in{\mathcal B}$.
\par
If ${\mathcal A},{\mathcal B}\subseteq \mc{P}(\mc{P}(X))$, then the {\emm covering property} U${}_{\scriptstyle fin}({\mathcal A},{\mathcal B})$ means: for any sequence of sets $\seq{{\mathcal U}_n\in{\mathcal A}}{n}$, for every $n$ there exists a~finite set ${\mathcal V}_n\subseteq {\mathcal U}_n$ such that $\{\bigcup{\mathcal V}_n:n\in\omega\}\in{\mathcal B}$ or there exists $n$ such that $X=\bigcup{\mathcal V}_n$.
\par
If ${\mathcal A},{\mathcal B}\subseteq \mc{P}(\mc{P}(X))$ and S$_1({\mathcal A},{\mathcal B})$ holds true, then we say that $X$ {\emm has the covering property} S$_1({\mathcal A},{\mathcal B})$. Sometimes a~topological space possessing the property S$_1({\mathcal A},{\mathcal B})$ is called S$_1({\mathcal A},{\mathcal B})$-space. Similarly for S${}_{\scriptstyle fin}({\mathcal A},{\mathcal B})$ and U$_{\scriptstyle fin}({\mathcal A},{\mathcal B})$.
\par
Let $\mc{F},\mc{G}\subseteq \mc{P}({}^X\RR)$ for some $X$. Then we prefer to speak about the~{\emm selection principle} S$_1({\mathcal F},{\mathcal G})$ or S${}_{\scriptstyle fin}({\mathcal F},{\mathcal G})$. Instead of U${}_{\scriptstyle fin}$ we define the {\emm selection principle} U${}_{\scriptstyle fin}^*(\mc{F},\mc{G})$ as follows: for any sequence of sets $\seq{F_n\in \mc{F}}{n}$, for every $n$ there exists a~finite set $H_n\subseteq F_n$ such that $\{\min H_n:n\in\omega\}\in\mc{G}$.
\par
If $F\subseteq {}^X\RR$ is infinite and $\Phi$ is one of the symbols ${\mathcal O}$, $\Omega$ or $\Gamma$, then we set
\[\Phi_{\snula}(F)=\{H\subseteq F: H\mbox{\ has the property\ }(\Phi_{\snula})\}.\]
If S${}_1(\Phi_{\snula}(F), \Psi_{\snula}(F))$ holds true, we shall say that $F$ {\emm satisfies the~selection principle} S${}_1(\Phi_{\snula}, \Psi_{\snula})$. We can identify the countable sets of functions with sequences of functions. If~S${}_1(\Phi_{\snula}(F)_{\scriptstyle cnt},\Psi_{\snula}(F)_{\scriptstyle cnt})$ holds true, we say that $F$ {\emm satisfies the~sequence selection principle} S${}_1(\Phi_{\snula},\Psi_{\snula})$. Similarly for the~selection principles S${}_{\scriptstyle fin}(\Phi_{\snula}, \Psi_{\snula})$ and U${}^*_{\scriptstyle fin}(\Phi_{\snula}, \Psi_{\snula})$.
\par
Note the following. By a~sequence selection principle we chose a set, not a sequence. In a~sequence a~function may occur infinitely many times. 
\par
Some known results about equivalences of introduced notions, which were  proved mainly in \cite{Sch} and \cite{JMSS},  may be easily generalized.
\begin{lemma}\label{OLambda}
For any family ${\mathcal A}$, for $\Phi=\Omega_{\scriptstyle cnt},\Gamma$ we have
\begin{eqnarray*}
{\rm S}_1(\Phi({\mathcal A}),{\mathcal O}({\mathcal A}))&\equiv&
{\rm S}_1(\Phi({\mathcal A}),\Lambda({\mathcal A})),\\  
{\rm S}_{\scriptstyle fin}(\Phi({\mathcal A}),{\mathcal O}({\mathcal A}))&\equiv& {\rm S}_{\scriptstyle fin}(\Phi({\mathcal A}),\Lambda({\mathcal A})),\\
{\rm U}_{\scriptstyle fin}(\Phi({\mathcal A}),{\mathcal O}({\mathcal A}))&\equiv&
{\rm U}_{\scriptstyle fin}(\Phi({\mathcal A}),\Lambda({\mathcal A})).
\end{eqnarray*}
\end{lemma}
\pf
We show the implication from left to right of the third line. Let a~sequence of $\varphi$-covers, all subsets of ${\mathcal A}$, be one-to-one enumerated as $\seq{{\mathcal U}_{n,m}}{n,m}$. If there exists a~finite ${\mathcal V}\subseteq{\mathcal U}_{n,m}$ such that $\bigcup{\mathcal V}=X$ then we have immediately ${\rm U}_{\scriptstyle fin}(\Phi,\Lambda)$. If $\Phi$ is $\Omega_{\scriptstyle cnt}$, then ${\mathcal U}_{n,m}=\{U^k_{n,m}:k\in\omega\}$.
If $\Phi=\Gamma$, then we take $\gamma$-subcovers $\{U^k_{n,m}:k\in\omega\}\subseteq {\mathcal U}_{n,m}$. In both cases we assume that the enumeration is one-to-one. 
\par
Throwing out from ${\mathcal U}_{n,m}$ the finite set
\[
\{U\in{\mathcal U}_{n,m}:(\exists i,j,l<n)\,U=U^l_{i,j}\},
\]
we obtain a~$\varphi$-cover ${\mathcal V}_{n,m}$. Now apply  
U${}_{\scriptstyle fin}(\Phi,{\mathcal O})$ to each sequence
$\langle {\mathcal V}_{n,m}:m\in\omega\rangle$. We obtain a~cover~${\mathcal W}_n$. By assumption $X\notin {\mathcal W}_n$. We claim that ${\mathcal W}=\bigcup_n{\mathcal W}_n$ is a~$\lambda$-cover.
Indeed, if $x\in Y\in {\mathcal W}_n$, then there exist $j,k$ such that $Y=U^k_{n,j}$. Let $n_0>n,j,k$. Then there exists a~set $Z\in {\mathcal W}_{n_0}$ such that $x\in Z$. By definition the set $Y$ was thrown out from every ${\mathcal U}_{n_0,m}$, so $Z\neq Y$.
\qed
\par
We shall need an~auxiliary result similar to Lemma~3.1 of \cite{BL2}.
\begin{lemma}\label{Setsversussequence}
Assume that the family ${\mathcal A}\subseteq {\mathcal P}(X)$ is closed under finite intersections. If $\varphi=\omega,\,\gamma$ then
\begin{enumerate}
\item[{\rm a)}] 
 {\rm S}${}_1(\Phi({\mathcal A}),\Gamma({\mathcal A}))$ is equivalent to the following:
\newline
for any sequence $\langle {\mathcal U}_n:n\in\omega\rangle$ of $\varphi$-covers, each being subset of~${\mathcal A}$, for every $n$ there exists a set $U_n\in {\mathcal U}_n$ such that
\begin{equation}\label{versus1}
(\forall x\in X)(\exists n_0)(\forall n\geq n_0)\, x\in U_n.
\end{equation}
\item[{\rm b)}] 
 {\rm S}${}_{\scriptstyle fin}(\Phi({\mathcal A}),\Gamma({\mathcal A}))$ is equivalent to the following:
\newline
for any sequence $\langle {\mathcal U}_n:n\in\omega\rangle$ of $\varphi$-covers, each being subset of~${\mathcal A}$, for every $n$ there exists a finite set ${\mathcal V}_n\subseteq {\mathcal U}_n$ such that
\begin{equation}\label{versus2}
(\forall x\in X)(\exists n_0)(\forall n\geq n_0)(\forall U\in {\mathcal V}_n)\, x\in U.
\end{equation}
\item[{\rm c)}] 
{\rm U}${}_{\scriptstyle fin}(\Phi({\mathcal A}),\Gamma({\mathcal A}))$ is equivalent to the following:
\newline
for any sequence $\langle {\mathcal U}_n:n\in\omega\rangle$ of $\varphi$-covers, each being subset of~${\mathcal A}$, for every $n$ there exists a finite set ${\mathcal V}_n\subseteq {\mathcal U}_n$ such that
\begin{equation}\label{versus3}
(\forall x\in X)(\exists n_0)(\forall n\geq n_0)\,x\in \bigcup{\mathcal V}_n
\end{equation}
or there exists an~$n$ such that $\bigcup {\mathcal V}_n=X$.
\end{enumerate}
\end{lemma}
\pf
\par
We show just the implications from left to right. The opposite implications are trivial.
\par
If $\varphi=\omega$ take an~$\omega$-cover ${\mathcal W}_n$ that is a~common refinement of ${\mathcal W}_0,\dots,{\mathcal W}_{n-1}$ and ${\mathcal U}_n$. If $\varphi=\gamma$ let ${\mathcal W}_0$ be a~countable $\gamma$-subcover of~${\mathcal U}_0$. For $n>0$  take a~countable $\gamma$-cover ${\mathcal W}_n$  which is a~common refinement of ${\mathcal W}_0,\dots,{\mathcal W}_{n-1}$ and ${\mathcal U}_n$. Since by \eqref{common}, two $\varphi$-covers, subsets of~${\mathcal A}$, have a~common refinement, again a~subset of ${\mathcal A}$, you can keep the covers ${\mathcal W}_n$ to be subsets of ${\mathcal A}$. 
\par
a) Apply S${}_1(\Phi({\mathcal A}),\Gamma({\mathcal A}))$ to the sequence $\langle {\mathcal W}_n:n\in\omega\rangle$. You obtain a~sequence $\seq{W_n}{n}\in\Gamma({\mathcal A})$, $W_n\in {\mathcal W}_n$ for each $n$. Then there exists an~increasing sequence of integers $\seq{n_k}{k}$ such that 
\[\{W_{n_k}:k\in\omega\}=\{W_n:n\in\omega\}\setminus \{\emptyset\}\]
and $W_{n_k}\neq W_{n_l}$ for $k\neq l$. For every $n_k\leq n<n_{k+1}$ take a~set $U_n\in{\mathcal U}_n$ such that $W_{n_{k+1}}\subseteq U_n$. If $n_0>0$ then take $U_n\in{\mathcal U}_n$ similarly for $n<n_0$. Since $\{W_{n_k}:k\in\omega\}$ is a~one-to-one enumeration of a~$\gamma$-cover, for every $x\in X$ there exists a~$k_0$ such that $x\in W_{n_{k+1}}$ for each $k\geq  k_0$. Then $x\in U_n$ for each $n\geq n_{k_0}$. 
\par
b) Apply S${}_{\scriptstyle fin}(\Phi({\mathcal A}),\Gamma({\mathcal A}))$ to the sequence $\langle {\mathcal W}_n:n\in\omega\rangle$. We obtain finite sets ${\mathcal Y}_n\subseteq {\mathcal W}_n$ such that $\bigcup_{n\in\omega} {\mathcal Y}_n$ is a~$\gamma$-cover. We construct a~one-to-one enumeration of $\bigcup_{n\in\omega} {\mathcal Y}_n$ as follows.
\par
Let $n_0$ be such that ${\mathcal Y}_0=\{Y_i:i<n_0\}$. By induction, $n_{k+1}$ is such that ${\mathcal Y}_{k+1}\setminus \bigcup_{l\leq k}{\mathcal Y}_l=\{Y_i:n_k\leq i< n_{k+1}\}$. The sequence $\seq{n_k}{k}$ is non-decreasing and  $\bigcup_{n\in\omega} {\mathcal Y}_n=\{Y_i:i\in\omega\}$.
\par
Since $Y_i\in {\mathcal Y}_{k}$ for $n_{k-1}\leq i< n_k$, there exists $U_i\in{\mathcal U}_k$ such that $Y_i\subseteq U_i$. Let ${\mathcal V}_k=\{U_i:n_{k-1}\leq i<n_k\}$. 
\par
Since $\{Y_i:i\in\omega\}$ is a~one-to-one enumeration of a~$\gamma$-cover, for every $x\in X$ there exists an~$i_0$ such that $x\in Y_i$ for $i\geq i_0$. Then, if $n_k\geq i_0$, we have $x\in U$ for each $U\in {\mathcal V}_k$.   
\par
c) Apply U${}_{\scriptstyle fin}(\Phi({\mathcal A}),\Gamma({\mathcal A}))$ to the sequence $\langle {\mathcal W}_n:n\in\omega\rangle$. We obtain finite sets ${\mathcal Y}_n\subseteq {\mathcal W}_n$ such that $\{\bigcup {\mathcal Y}_n:n\in\omega\}$ is a~$\gamma$-cover or 
$\bigcup {\mathcal Y}_n= X$ for some $n$. We may assume that $\bigcup {\mathcal Y}_n\neq X$ for each $n$. 
Then there exists an~increasing sequence of integers $\seq{n_k}{k}$ such that $\{\bigcup {\mathcal Y}_{n_k}:k\in\omega\}=\{\bigcup {\mathcal Y}_n:n\in\omega\}\setminus \{\emptyset\}$ and $\bigcup {\mathcal Y}_{n_k}\neq \bigcup {\mathcal Y}_{n_l}$ for $k\neq l$. For every $n_k\leq n<n_{k+1}$ take finite set ${\mathcal V}_n\subseteq {\mathcal U}_n$ such that each $Y\in {\mathcal Y}_{n_{k+1}}$ is subset of some $V\in {\mathcal V}_n$.
\par
Since $\{\bigcup {\mathcal Y}_{n_k}:k\in\omega\}$ is a~one-to-one enumeration of a~$\gamma$-cover, for every $x\in X$ there exists a~$k_0$ such that $x\in \bigcup {\mathcal Y}_{n_k}$ for each $k>k_0$. Then $x\in \bigcup {\mathcal V}_n$ for each $n\geq n_{k_0}$. 
\qed  
\par
In~\cite{JMSS} the authors show that several S${}_{\scriptstyle fin}$ properties are equivalent to some of S${}_1$ or U${}_{\scriptstyle fin}$ properties, at least for countable covers.
For example, using Lemma~\ref{Setsversussequence} one can easily see that for any family ${\mathcal A}\subseteq {\mathcal P}(X)$ closed under finite intersections we have
\[{\rm S}_{\scriptstyle fin}(\Phi({\mathcal A}),\Gamma({\mathcal A}))={\rm S}_1(\Phi({\mathcal A}),\Gamma({\mathcal A}))\textnormal{\ for\ }\Phi=\Omega,\,\Gamma.\]
Indeed, let $\seq{{\mathcal U}_n}{n}$ be a~sequence of $\varphi$-covers such that ${\mathcal U}_n\subseteq {\mathcal A}$ for each $n$. Since ${\mathcal A}$ is closed under finite intersections, for every $n$ there exists a~$\varphi$-cover ${\mathcal V}_n$ such that ${\mathcal V}_n$ is a~refinement of ${\mathcal U}_n$ and, if $n>0$, also of ${\mathcal V}_{n-1}$. Applying ${\rm S}_{\scriptstyle fin}(\Phi({\mathcal A}),\Gamma({\mathcal A}))$ to the sequence $\seq{{\mathcal V}_n}{n}$, for every $n$ we obtain finite set ${\mathcal W}_n\subseteq {\mathcal V}_n$ such that $\bigcup_n{\mathcal W}_n$ is a~$\gamma$-cover. Let $\seq{n_k}{k}$ be the increasing sequence of those $n$ for which ${\mathcal W}_n\neq\emptyset$. For every $k$, choose a~$W_k\in{\mathcal W}_{n_k}$. For $n_k\leq n<n_{n_{k+1}}$,
choose a~$U_n\in{\mathcal U}_n$ such that $W_{k+1}\subseteq U_n$. Then $\{U_n:n\in \omega\}$ is the desired $\gamma$-cover.  
\par  
If ${\mathcal A}$ is closed under finite unions, then we can show that 
\[{\rm S}_{\scriptstyle fin}(\Gamma({\mathcal A}),{\mathcal O}({\mathcal A}))={\rm U}_{\scriptstyle fin}(\Gamma({\mathcal A}),{\mathcal O}({\mathcal A})).\] 
\par
Let $\Psi$ denote one of the symbols ${\mathcal O}$, $\Lambda$, $\Omega$ or $\Gamma$, and let ${\mathcal A}$ be closed under finite unions. Then one can easily see that
\begin{eqnarray}\label{Ufineq}
{\rm U}_{\scriptstyle fin}({\mathcal O}({\mathcal A}),\Psi({\mathcal A}))&\equiv&
{\rm U}_{\scriptstyle fin}(\Omega({\mathcal A}),\Psi({\mathcal A})),\\
{\rm U}_{\scriptstyle fin}({\mathcal O}({\mathcal A})_{\scriptstyle cnt},\Psi({\mathcal A}))&\equiv&{\rm U}_{\scriptstyle fin}(\Gamma({\mathcal A}),\Psi({\mathcal A})).\nonumber
\end{eqnarray}
However for uncountable covers
\begin{equation}\label{nonUnifeq}
{\rm U}_{\scriptstyle fin}(\Gamma({\mathcal T}),\Psi({\mathcal T}))\to {\rm U}_{\scriptstyle fin}(\Omega({\mathcal T}),\Psi({\mathcal T}))\ does\ not\ hold\ true.
\end{equation}
Indeed, a~countably compact non-compact topological space has the property U${}_{\scriptstyle fin}(\Gamma,\Gamma)$ and does not
have the property U${}_{\scriptstyle fin}({\mathcal O},{\mathcal O})$.
$\beta\omega\setminus\{x\}$, where $x\in \beta\omega\setminus \omega$, is a~countably compact non-compact topological space -- see \cite{Eng}, p. 262. The existence of a~perfectly normal countably compact non-compact topological space is undecidable in {\bf ZFC} -- see \cite{Vau}, pp. 597--598.
%%%%%%%%%%%%%%%%%%%%%%%%%%%%%%%%%%%%%%%%%%%%%%%%%%%%%%%%%%%%%%
%%%%%%%%%%%%%%%%%%%%%%%%%%%%%%%%%%%%%%%%%%%%%%%%%%%%%%%%%%%%%%
\section{Covering Properties versus Selection Principles}\label{CPCP}
In~\cite{BL2} we have shown connections (in some cases already known) of the property S${}_1$ of open covers or open shrinkable covers of a~topological space~$X$ with the corresponding selection principles for $\CP{X}$ or $\USC{X}^+$, respectively. We can easily prove similar results for the S${}_{\scriptstyle fin}$ and U${}_{\scriptstyle fin}$ properties.
\par
Now we extend those results for covers from a~family ${\mathcal A}$ of sets.
\begin{theorem}\label{Ufinsemi}
Assume that $\Phi$ is one of the symbols  $\Omega$ or $\Gamma$ and $\Psi$ is one of the symbols ${\mathcal O}$, $\Omega$ or $\Gamma$. Let ${\mathcal A}$ be a~family of subsets of a~set~$X$ closed under finite unions. If $\Psi=\Gamma$, we assume that ${\mathcal A}$ is also closed under finite intersections. Then for any couple $\langle\Phi,\Psi\rangle$ different from $\langle\Omega,{\mathcal O}\rangle$ a~topological space $X$ possesses the covering property~{\rm S}${}_1(\Phi({\mathcal A}),\Psi({\mathcal A}))$ if and only if $\USM{X}{{\mathcal A}}^+$ satisfies the selection principle~{\rm S}${}_1(\Phi_{\snula},\Psi_{\snula})$.
\par
If ${\mathcal A}$ is a~$\sigma$-algebra of sets, then the family $\USM{X}{{\mathcal A}}^+$ may be replaced by $\Meas{X}{{\mathcal A}}$. 
\end{theorem} 
\pf
\par
We assume that $X$ has the covering property {\rm S}${}_1(\Phi({\mathcal A}),\Psi({\mathcal A}))$, $\Phi=\Omega,\Gamma$ and $\Psi={\mathcal O},\Omega,\Gamma$, $\langle\Phi,\Psi\rangle\neq\langle\Omega,{\mathcal O}\rangle$. Let $\seq{F_n}{n}$ be a~sequence of infinite sets of non-negative upper ${\mathcal A}$-semimeasurable real functions with the property $(\Phi_{\snula})$. By the assumption of upper semimeasurability we obtain that $\mc{U}^{\varepsilon}(F_n)\subseteq {\mathcal A}$ for any $n$ and any $\varepsilon>0$. 
\par
If  $\Phi=\Omega$, then by Theorem~\ref{gamma} for every $n$ either there exists a~countable subset of $F_n$ uniformly converging to $\nula$ or there exists a~$\delta_n>0$ such that $\mc{U}^{\varepsilon}(F_n)$ is an~$\omega$-cover for every positive $\varepsilon<\delta_n$.
\par
If $\Phi=\Gamma$, then by Theorem~\ref{gamma} for every $n$ either $F_n$ is countable and $F_n\rightrightarrows \nula$ or there exists $\delta_n>0$ such that $\mc{U}^{\varepsilon}(F_n)$ is a~$\gamma$-cover for each positive $\varepsilon<\delta_n$ and for every $f\in F_n$ the set $\{g\in F_n:U^{\varepsilon}_g=U^{\varepsilon}_f\}$ is finite.
\par
In both cases, $\Phi=\Omega$ or $\Phi=\Gamma$, let $A$ be the set of those $n\in\omega$ for which there exists a~$\delta_n>0$ such that for every $\varepsilon<\delta_n$ the family ${\mathcal U}^{\varepsilon}(F_n)$ is a~$\varphi$-cover.
\par
If $A$ is finite, one can find a~sequence $\langle f_n\in F_n:n\in \omega\setminus A\rangle$ such that $f_n\rightrightarrows\nula$. The family $\{f_n:n\in\omega\setminus A\}$ has the property~$(\Gamma_{\snula})$. For $n\in A$ take $f_n\in F_n$, $f_n\neq\nula$, arbitrary. Then the family $\fml{f_n}{n}$ has the property $(\Gamma_{\snula})$ as well. If $\Psi=\Omega$ then $f_n\neq \nula$ for every $n$ and hence the family $\fml{f_n}{n}$ has the property $(\Omega_{\snula})$.
\par
So, let $A$ be infinite. We set $\varepsilon_n=\min\{\delta_n/2,2^{-n}\} \mbox{\rm\ for\ }n\in A$.
\par
If $\Psi=\Omega$, we apply S${}_1(\Phi,\Omega)$ to the sequence \mbox{$\{{\mathcal U}^{\varepsilon_n}(F_n): n\in A\}$}. We obtain sets $U_n\in{\mathcal U}^{\varepsilon_n}(F_n)$, $n\in A$ such that $\{U_n: n\in A\}$ is an~$\omega$-cover. For every $n\in A$ there exists a~function $f_n\in F_n$ such that $U_n=U^{\varepsilon_n}_{f_n}$. By Lemma~\ref{epsilontozero} either there exists a~subsequence of $\{f_n:n\in A\}$  uniformly converging to~$\nula$ or there exists a~$\delta>0$ such that the family $\{U_{f_n}^{\varepsilon}:n\in A\}$ is an~$\omega$-cover for each positive $\varepsilon<\delta$. Therefore the family $\{f_n:n\in A\}$ has the property~$(\Omega_{\snula})$.
\par
Now, let $\Psi=\Gamma$. Since for every $n\in A$ the family $\{{\mathcal U}^{\varepsilon_n}(F_n)$ is \mbox{a~$\gamma$-cover,} by Lemma~\ref{Setsversussequence} for every $n\in A$ there exists~$U_n\in {\mathcal U}^{\varepsilon_n}(F_n)$ such that
\[(\forall x\in X)(\exists n_0)(\forall n\geq n_0,n\in A)\, x\in U_n.\]
For every $n\in A$ there exists a~function $f_n\in F_n$ such that $U_n=U^{\varepsilon_n}_{f_n}$. Since $\varepsilon_n\to 0$, we obtain that there exists an~integer~$n_0$ such that $\varepsilon_n\leq \varepsilon$ and $x\in U^{\varepsilon_n}_{f_n}$ for each $n\geq n_0,n\in A$. Then $f_n(x)<\varepsilon$ for each $n\geq n_0,n\in A$. Thus  the family $\{f_n:n\in A\}$ possesses the property $(\Gamma_{\snula})$.  
\par
If $\Psi={\mathcal O}$, then $\Phi=\Gamma$ and by Lemma~\ref{OLambda} S${}_1(\Gamma,{\mathcal O})$ is equivalent to ${\rm S}{}_1(\Gamma,\Lambda)$. Applying  S${}_1(\Gamma,\Lambda)$ to the sequence $\{{\mathcal U}^{\varepsilon_n}(F_n): n\in A\}$ we obtain sets $U_n\in{\mathcal U}^{\varepsilon_n}(F_n)$, $n\in A$ such that the family $\{U_n: n\in A\}$ is \mbox{a~$\lambda$-cover.} We apply Lemma~\ref{epsilontozero} to the family $\{U_n: n\in A\}$ and continue as above.  
\par
If $\Psi$ is ${\mathcal O}$ or $\Omega$, for $n\notin A$ take $f_n\in F_n$ arbitrary. Then the family $\fml{f_n}{n}$ has the property $(\Psi_{\snula})$ as well.
\par
Assume that $\Psi=\Gamma$. If $\omega\setminus A$ is finite, similarly as above, take $f_n\in F_n$ for $n\in\omega\setminus A$ arbitrary. Then $\fml{f_n}{n}$ has the property $(\Gamma_{\snula})$. So let $\omega\setminus A$ be infinite. One can easily construct a~sequence $\langle f_n\in F_n:n\in\omega\setminus A\rangle$ such that $f_n\rightrightarrows 0$, $n\in\omega\setminus A$. Then the family $\{f_n:n\in\omega\}$ has the property $(\Gamma_{\snula})$. 
%%%%%%%%%%%%%%%%%%%%%%%%%%%%%%%%%%%%%%%%%%%%%
\par
Let $\USM{X}{{\mathcal A}}^+$ satisfy the selection principle~{\rm S}${}_1(\Phi_{\snula},\Psi_{\snula})$. Assume that $\seq{\mc{U}_n }{n}$ is a~sequence of $\varphi$-covers, ${\mc{U}_n }\subseteq {\mathcal A}$ for each $n\in\omega$.
For any~set $U\in{\mathcal A}$ we set
\begin{equation}\label{fuscnm}
f_U(x)=\left\{
\begin{array}{ll}
0&\mbox{\ if\ } x\in U,\\
1&\mbox{\ otherwise}.
\end{array}
\right.
\end{equation}
Then every $f_U$ is a~non-negative upper ${\mathcal A}$-semicontinuous real functions.
If ${\mathcal A}$ is an~algebra of sets, then $f_U$ is ${\mathcal A}$-measurable for each $U\in{\mathcal A}$.
\par
We set $F_n=\{f_U:U\in\mc{U}_n\}$. Then $F_n$ is a~family of non-negative upper ${\mathcal A}$-semicontinuous real functions.
\par
For every $\varepsilon\leq 1$, every $n\in\omega$, every $U\in {\mathcal U}_n$ we have $U_{f_U}^{\varepsilon}=U$. Thus for every $\varepsilon\leq 1$ every ${\mathcal U}^{\varepsilon}(F_n)$ is a~$\varphi$-cover.
\par
By Theorem~\ref{gamma}, every set $F_n$ has the property ($\Phi_{\snula}$). By  the selection principle U${}_{\scriptstyle fin}^*(\Phi_{\snula}({\rm USC}^+),\Psi_{\snula}({\rm USC}^+))$ for every $n\in\omega$ 
there exists $f_{n,0},\dots,f_{n,k_n}\in F_n$ such that the set 
\[\fml{\min\{ f_{n,0},\dots,f_{n,k_n}\}}{n}\]
has the property~($\Psi_{\snula}$). For every $n$ and $i\leq k_n$ there exist a~set $U_{n,i}\in {\mathcal U}_n$ such that $f_{n,i}=f_{U_{n,i}}$. By Theorem~\ref{gamma}, either there exists a~subsequence of the family $\fml{\min\{ f_{n,0},\dots,f_{n,k_n}\}}{n}$ uniformly converging to $\nula$ -- then $\bigcup_{i\leq k_n}U_{n,i}=X$ for some $n$ -- or there exists a~$\delta$ such that for every $\varepsilon<\delta$ the family 
\[{\mathcal U}^{\varepsilon}(\fml{\min\{ f_{n,0},\dots,f_{n,k_n}\}}{n})\]
is a~$\psi$-cover. Since $U^{\varepsilon}_{\min\{f_{n,0},\dots,f_{n,k_n}\}}=\bigcup_{i\leq k_n}U_{n,i}$ for $\varepsilon<1$, the family $\{\bigcup_{i\leq k_n}U_{n,i}:n\in\omega\}$ is a~$\psi$-cover as well.  
\qed
%%%%%%%%%%%%%%%%%%%%%%%%%%%%%%%%%%%%%%%%%%%%%%%%%%
%%%%%%%%%%%%%%%%%%%%%%%%%%%%%%%%%%%%%%%%%%%%%%%%%%
\par
Modifying the proof of Theorem~\ref{Ufinsemi} similarly as in \cite{BL2} one obtains
\begin{theorem}\label{UfinsemiCont}
Assume that $\Phi$ is one of the symbols $\Omega$ or $\Gamma$and $\Psi$ is one of the symbols ${\mathcal O}$, $\Omega$ or $\Gamma$. Then for any couple $\langle\Phi,\Psi\rangle$ different from $\langle\Omega,{\mathcal O}\rangle$ a~normal topological space~$X$ has the property~{\rm S}${}_1(\Phi^{\scriptstyle sh},\Psi)$, {\rm S}${}_{\scriptstyle fin}(\Phi^{\scriptstyle sh},\Psi)$, {\rm U}${}_{\scriptstyle fin}(\Phi^{\scriptstyle sh},\Psi)$ respectively, if and only if $\CP{X}$ satisfies the selection principle~{\rm S}$_1(\Phi_{\snula},\Psi_{\snula})$, {\rm S}${}_{\scriptstyle fin}(\Phi_{\snula},\Psi_{\snula})$, {\rm U}${}_{\scriptstyle fin}^*(\Phi_{\snula},\Psi_{\snula})$, respectively.
\end{theorem} 
\pf
\par
If $F$ is a~family of continuous functions with property $(\Phi_{\snula})$  then the family $\vert F\vert$ is a~family of non-negative continuous functions and by \eqref{absolute}, possesses the property $(\Phi_{\snula})$ as well. So the implication from left to right follows by Theorem~\ref{Ufinsemi} for ${\mathcal A}={\mathcal T}$.
\par  
If $\seq{\mathcal V_n}{n}$, $\seq{\mathcal U_n}{n}$ are sequences of open~$\varphi$-covers such that
\[(\forall V\in {\mathcal V_n})(\exists U_V\in {\mathcal V_n})\, \overline{V}\subseteq U_V,\] 
then for any $V\in{\mathcal V}_n$ there exists a~continuous function $f_v$ such that 
\begin{equation*}
f_V(x)=\left\{
\begin{array}{ll}
0&\mbox{\ if\ } x\in V,\\
1&\mbox{\ if\ } x\in X\setminus U_V.
\end{array}
\right.
\end{equation*}
Apply the selection principle~{\rm S}$_1(\Phi_{\snula},\Psi_{\snula})$ to the sequence $\seq{F_n}{n}$, where
\[F_n=\{f_V:V\in {\mathcal V}_n\}.\]
Similarly for {\rm S}${}_{\scriptstyle fin}(\Phi_{\snula},\Psi_{\snula})$ and {\rm U}${}_{\scriptstyle fin}^*(\Phi_{\snula},\Psi_{\snula})$.
\qed
\par
By a~modification of the proof of Theorem~7.4 of \cite{BL2} we obtain
\begin{theorem}\label{omegaequiv}
Let $X$ be a~Tychonoff topological space,  $\Psi$ being one of the symbols ${\mathcal O}, \Omega, \Gamma$. Then the following are equivalent:
\begin{enumerate}
\item[{\rm (i)}] $X$ is an~{\rm S}${}_{\scriptstyle fin}(\Omega,\Psi)$-space.
\item[{\rm (ii)}] $\USC{X}^+$ satisfies the selection principle {\rm S}${}_{\scriptstyle fin}(\Omega_{\snula},\Psi_{\snula})$.
\item[{\rm (iii)}] $\CP{X}$ satisfies the selection principle {\rm S}${}_{\scriptstyle fin}(\Omega_{\snula},\Psi_{\snula})$.
\end{enumerate}
Similarly for {\rm U}${}_{\scriptstyle fin}(\Omega,\Psi)$.
\end{theorem}
\par
We can prove accordingly modified Theorems for countable covers and sequence selection principles. Moreover, in the case of countable covers we do not need to omit the couple $\langle \Omega,{\mathcal O}\rangle$, since we can use the first equivalence of~Lemma~\ref{OLambda}.
\begin{theorem}\label{Ufinsemicount}
Assume that $\Phi$ is one of the symbols $\Omega$ or $\Gamma$ and $\Psi$ is one of the symbols ${\mathcal O}$, $\Omega$ or $\Gamma$.  Let ${\mathcal A}$ be a~family of subsets of a~set~$X$ closed under finite unions. If $\Psi=\Gamma$ we assume also that ${\mathcal A}$ is closed under finite intersections. Then for any couple $\langle\Phi,\Psi\rangle$ the following holds true:
\begin{enumerate}
\item[{\rm a)}]
$X$is has the covering property~{\rm S}${}_1(\Phi({\mathcal A})_{\scriptstyle cnt},\Psi({\mathcal A}))$ if and only if the set $\USM{X}{{\mathcal A}}^+$ satisfies the sequence selection principle~{\rm S}${}_1(\Phi_{\snula},\Psi_{\snula})$.
\item[{\rm b)}]
$X$ has the covering property~{\rm S}${}_{\scriptstyle fin}(\Phi({\mathcal A})_{\scriptstyle cnt},\Psi({\mathcal A}))$ if and only if the set $\USM{X}{{\mathcal A}}^+$ satisfies the sequence selection principle~{\rm S}${}_{\scriptstyle fin}(\Phi_{\snula},\Psi_{\snula})$.
\item[{\rm c)}]
$X$ has the covering property~{\rm U}${}_{\scriptstyle fin}(\Phi({\mathcal A})_{\scriptstyle cnt},\Psi({\mathcal A}))$ if and only if the set $\USM{X}{{\mathcal A}}^+$ satisfies the sequence selection principle~{\rm U}${}_{\scriptstyle fin}^*(\Phi_{\snula},\Psi_{\snula})$.
\end{enumerate}
If ${\mathcal A}$ is a~$\sigma$ algebra of sets, then the family $\USM{X}{{\mathcal A}}^+$ may be replaced by $\Meas{X}{{\mathcal A}}$. 
\end{theorem} 
\par
Let {\Borel}($X$), or simply {\Borel}, be the family of all Borel subsets of the~topological space $X$. We obtain
\begin{corollary}
Assume that $\Phi$ is one of the symbols $\Omega$, $\Gamma$ and $\Psi$ is some of the symbols ${\mathcal O}$, $\Omega$, $\Gamma$. Then for any couple $\langle\Phi,\Psi\rangle$ different from $\langle\Omega,{\mathcal O}\rangle$ the~topological space $X$ has the covering property {\rm S}$_1(\Phi(\mbox{\Borel}),\Psi(\mbox{\Borel}))$ if and only if the set of all Borel measurable real functions satisfies the selection principle {\rm S}${}_1(\Phi_{\snula},\Psi_{\snula})$.
\end{corollary}
\par
For the covering property {\rm S}$_{\scriptstyle fin}(\Phi(\mbox{\Borel}),\Psi(\mbox{\Borel}))$, one can prove similar results. We have also
\begin{corollary}\label{UfinBorel}
Assume that $\Phi$ is one of the symbols $\Omega$, $\Gamma$ and $\Psi$ is some of the symbols ${\mathcal O}$, $\Omega$, $\Gamma$. Then for any couple $\langle\Phi,\Psi\rangle$ different from $\langle\Omega,{\mathcal O}\rangle$ a~topological space $X$ has the covering property {\rm U}$_{\scriptstyle fin}(\Phi(\mbox{\Borel}),\Psi(\mbox{\Borel}))$ if and only if the set of all Borel measurable real functions satisfies the selection principle {\rm U}${}_{\scriptstyle fin}^*(\Phi_{\snula},\Psi_{\snula})$.
\end{corollary}
%%%%%%%%%%%%%%%%%%%%%%%%%%%%%%%%%%%%%%%%%%%%%%%%%%%%%%%%%%%%%%%%
%%%%%%%%%%%%%%%%%%%%%%%%%%%%%%%%%%%%%%%%%%%%%%%%%%%%%%%%%%%%%%%%
\section{Hurewicz like Theorem}
\par
We shall need a~result from~\cite{BH}.
\begin{lemma}[J. Hale\v s]\label{hales}
Let ${\mathcal A}$ be a~$\sigma$-algebra of subsets of $X$. If $\seq{f_n}{n}$ is a~sequence of positive ${\mathcal A}$-measurable functions converging monotonically to $\nula$, then there exists an~${\mathcal A}$-measurable  function \mbox{$h:X\longrightarrow \RR^+$} such that
\[(\forall x\in X)(\forall n>h(x))\,f_n(x)<1.\] 
\end{lemma}
\prfwp is easy. Set
\[g(x)=\sum_{n=0}^{\infty}2^{-n}\cdot\min\{1,f_n(x)\}\]
and $h(x)=\log_2(2-g(x))$.
\qed
\par
Now we can prove results similar to those of Hurewicz~\cite{Hu} and Tsaban~\cite{Ts} mentioned in Section~\ref{Int}.
\begin{theorem}\label{buksup}
Let ${\mathcal A}$ be a~$\sigma$-algebra of subsets of $X$ containing all open sets. Then the following are equivalent:
\begin{enumerate}
\item[{\rm (i)}] Every ${\mathcal A}$-measurable image of $X$ into ${}^{\omega}\RR$ is not dominating, not finitely dominating, eventually bounded, respectively.
\item[{\rm (ii)}] Every ${\mathcal A}$-measurable image of $X$ into ${}^{\omega}\omega$ is not dominating, not finitely dominating, eventually bounded, respectively.
\item[{\rm (iii)}] $X$ has the property {\rm U}${}_{\scriptstyle fin}(\Gamma({\mathcal A}),\Phi({\mathcal A}))$ for $\Phi={\mathcal O}$, $\Omega,\ \Gamma$, respectively.
\item[{\rm (iv)}] The set of all ${\mathcal A}$-measurable real functions  satisfies the selection principle {\rm U}${}_{\scriptstyle fin}^*(\Gamma_{\snula},\Phi_{\snula})$ for $\Phi={\mathcal O},\ \Omega,\ \Gamma$, respectively.
\end{enumerate}
Assuming that $X$ is zero dimensional, the equivalencies hold true also for ${\mathcal A}={\mathcal T}$, i.e., for continuous mappings.
\end{theorem}
\pf
The equivalence ${\rm (iii)}\equiv {\rm (iv)}$ is a~special case of Theorem~\ref{Ufinsemi}.
\par
Let $[x]$ denote the integer part of the real $x$. For $\alpha\in{}^{\omega}\RR$ we set 
\[p(\alpha)=\seq{\max\{0,[\alpha(n)]+1\}}{n}\in{}^{\omega}\omega.\]
The function $p:{}^{\omega}\RR\longrightarrow {}^{\omega}\omega$ is Borel measurable. If for $Y\subseteq {}^{\omega}\RR$ the image $p(Y)\subseteq {}^{\omega}\omega$ is eventually bounded by $\beta \in {}^{\omega}\omega$, then $Y$ is eventually bounded by same function $\beta\in{}^{\omega}\RR$. Similarly, if $p(Y)$ is not dominating or not finitely dominating, then $Y$ is also such. Since ${\mathcal A}$ is a~$\sigma$-algebra containing all open sets, we obtain that  if $f:X\longrightarrow {}^{\omega}\RR$ is ${\mathcal A}$-measurable, then $f^{-1}(B)\in {\mathcal A}$ for any Borel set $B\subseteq {}^{\omega}\RR$. Thus $f\circ p:X\longrightarrow {}^{\omega}\omega$ is ${\mathcal A}$-measurable.  Hence $({\rm ii})\to ({\rm i})$. The implication $({\rm i})\to ({\rm ii})$ is trivially true.
\par
We show ${\rm (i)}\equiv {\rm (iv)}$ for $\Phi=\Omega$. The other cases can be proved similarly. 
\par
Assume (i). Let $F_n=\seq{f_{n,m}}{m}$, $n\in\omega$ be a~sequence of sequences of ${\mathcal A}$-measurable functions such that $F_n\in\Gamma_{\snula}({\mathcal A})$ for every $n$, i.e., $f_{n,m}\to\nula$. We set $g_{n,m}=\min\{\vert f_{n,i}\vert :i\leq m\}$. Then the sequence $\seq{g_{n,m}}{m}$ converges monotonically to $\nula$ for each $n$.  
By Lemma~\ref{hales} for every $n$ there exists an~${\mathcal A}$-measurable function $h_n$ such that
\[(\forall x\in X)(\forall m\geq h_n(x))\,g_{n,m}(x)<2^{-n}.\] 
We set
\[h(x)=\seq{h_n(x)}{n}\in {}^{\omega}\RR\mbox{\ for\ }x\in X.\]
Thus $h_n(x)=h(x)(n)$. By (i) there exists a~function $\alpha\in {}^{\omega}\RR$ such that  
\[(\forall B\in[ X]^{<\omega})(\forall n_0)(\exists n\geq n_0)(\forall x\in B)\,h_n(x)\leq \alpha(n).\]
We set 
\[H_n=\{\vert f_{n,0}\vert,\dots,\vert f_{n,\alpha(n)}\vert\}.\]
Let $\varepsilon>0$ and $x_0,\dots,x_k\in X$. Then there exists sufficiently great $n$ such that $2^{-n}<\varepsilon$ and $h_n(x_i)\leq \alpha(n)$ for every $i=0,\dots,k$. Then $\min H_n(x_i)=g_{n,\alpha(n)}(x_i)<2^{-n}<\varepsilon$ for $i=0.,\dots,k$. Hence $\min H_n\in N^{\varepsilon}_{x_0,\dots,x_k}$.
\par
Assume (iv). Let $h:X\longrightarrow {}^{\omega}\RR$ be ${\mathcal A}$-measurable. We set $g_{n,m}(x)=\vert h(x)(n)\vert\cdot2^{-m}$. Then for every $n$, $\seq{g_{n,m}}{m}$ converges monotonically to $\nula$. 
Thus, by (iv) there exists a sequence $\seq{m_n}{n}$ such that
$\nula\in\overline{\{g_{n,m_n}:n\in\omega\}}$. Then for every~finite subset $\{x_0,\dots,x_k\}$ of $X$ there exists arbitrary large $n$ such that $g_{n,m_n}(x_i)<1$ for every $i=0,\dots,k$. Then $\vert h(x_i)(n)\vert <2^{m_n}$ for $i=0,\dots,k$. Hence
\[(\forall B\in[ X]^{<\omega})(\forall n_0)(\exists n\geq n_0)(\forall x\in B)\,h(x)(n)\leq 2^{m_n},\]
i.e., $h(X)$ is not finitely dominating.
\qed
\par
For the notion of a~QN-space and its properties see~\cite{BRR}.
\begin{corollary} For a~perfectly normal topological space $X$ the following are equivalent:
\begin{enumerate}
\item[{\rm (i)}] $X$ is a~{\rm QN}-space.
\item[{\rm (ii)}] Every Borel image of $X$ into ${}^{\omega}\RR$ is eventually bounded.  
\item[{\rm (iii)}] Every Borel image of $X$ into ${}^{\omega}\omega$ is eventually bounded.
\item[{\rm (iv)}] The set of all Borel measurable functions satisfies the selection principle {\rm U}${}_{\scriptstyle fin}^*(\Gamma_{\snula},\Gamma_{\snula})$.
\item[{\rm (v)}] $X$ is a~{\rm U}${}_{\scriptstyle fin}(\Gamma(\mbox{\Borel}),\Gamma(\mbox{\Borel}))$-space.
\end{enumerate}
\end{corollary}
\pf
The equivalencies $({\rm ii})\equiv ({\rm iii})\equiv ({\rm iv})\equiv ({\rm v})$ are special cases of Theorem~\ref{buksup}. The equivalence $({\rm i})\equiv ({\rm iii})$ is the Tsaban -- Zdomskyy's Theorem, see \cite{TZ}. 
\qed 
\par
Similarly as M.~Scheepers and B.~Tsaban~\cite{ST} we obtain
\begin{corollary}\label{SvchTs}
Let ${\mathcal A}$ be a~$\sigma$-algebra of subsets of $X$ containing all open sets. Then the following are equivalent:
\begin{enumerate}
\item[{\rm (i)}] $X$ has the property {\rm S}${}_1(\Gamma({\mathcal A}),\Phi({\mathcal A}))$ for $\Phi={\mathcal O}$, $\Omega,\ \Gamma$, respectively.
\item[{\rm (ii)}] $X$ has the property {\rm S}${}_{\scriptstyle fin}(\Gamma({\mathcal A}),\Phi({\mathcal A}))$ for $\Phi={\mathcal O}$, $\Omega,\ \Gamma$, respectively.
\item[{\rm (iii)}] $X$ has the property {\rm U}${}_{\scriptstyle fin}(\Gamma({\mathcal A}),\Phi({\mathcal A}))$ for $\Phi={\mathcal O}$, $\Omega,\ \Gamma$, respectively.
\item[{\rm (iv)}] Every ${\mathcal A}$-measurable image of $X$ into ${}^{\omega}\RR$ is not dominating, not finitely dominating, eventually bounded, respectively.
\end{enumerate}
\end{corollary}
\pf The implications ${\rm (i)}\to{\rm (ii)}\to {\rm (iii)}$ are trivial. Implication ${\rm (iii)}\to {\rm (iv)}$ is Theorem~\ref{buksup}. We  show the implication ${\rm (iv)}\to{\rm (i)}$ as in~\cite{ST}.
\par
 Let $\seq{\mc{U}_n }{n}$ be a~sequence of $\gamma$-covers,  for every $n\in\omega$ being ${\mc{U}_n }=\{U_n^m:m\in\omega\}\subseteq {\mathcal A}$. We define an~${\mathcal A}$-measurable function $f:X\longrightarrow {}^{\omega}\omega$ as follows:
\[f(x)(n)=\min\{k:(\forall m\geq k)\, x\in U^m_n\}.\]
If $\alpha\in{}^{\omega}\omega$ witnesses that $f(X)\subseteq {}^{\omega}\omega$ is not dominating, not finitely dominating, bounded, respectively, then $\{U_n^{\alpha(n)}:n\in\omega\}$ is a~$\varphi$-cover for $\Phi={\mathcal O}$, $\Omega,\ \Gamma$, respectively.  
\qed
%%%%%%%%%%%%%%%%%%%%%%%%%%%%%%%%%%%%%%%%%%%%%%%%%%%%%%%%%%%%%%%%%%%%%%%%
\section{Remarks}
 Several authors restrict their investigations to the countable covers. That is not our case. We consider covers which may be uncountable. \par
As we have already mentioned, many special cases of presented results are known. I will mention all those of them that I know. I am sorry if I have omitted some other known results. The only reason for such an~omitting is the fact that I do not know the result.
\par
For the covering property S${}_1$ we have presented a survey of recent results  in~\cite{BL2}.   
\par
For ${\mathcal A}={\mathcal T}$, $\Phi=\Gamma$ and $\Psi=\Omega$, the parts a)  and b) of Theorem~\ref{Ufinsemi} were proved by M.~Sakai~\cite{Sa}. The part a) of Theorem~\ref{Ufinsemi} for  ${\mathcal A}={\mathcal T}$ is Theorem~6.2 of~\cite{BL2}.
\par
The equivalence ${\rm (ii)}\equiv {\rm (iii)}$ of Theorem~\ref{buksup} for ${\mathcal A}=\Borel$ and $\Phi={\mathcal O},\,\,\Omega$ are Theorems~6. and 9. of M.~Scheepers and B.~Tsaban~\cite{ST}. For ${\mathcal A}=\Borel$ and $\Phi=\Gamma$, this equivalence is Theorem~3. of T.~Bartoszynski and M.~Scheepers~\cite{BS}. For ${\mathcal A}={\mathcal T}$ and $X=\RR$, L.S.~Zdomskyy~\cite{Zd} proved the equivalence $(ii)\equiv(iii)$ of Theorem~\ref{buksup}. 
\par
Finally, let me allow to mention about the terminology and mainly the notations which I have used in the paper. I suppose that the terminology is close to that used by lots of mathematicians working in this area.
\par
Some authors define a~cover ${\mathcal U}$ without asking $X\notin {\mathcal U}$. One can easily reformulate our results for such notion of a~cover. Maybe, that some results can be expressed in a~simpler way using such notion. For me, it seemed to be more convenient to assume that $X\notin {\mathcal U}$.  
\par
 I consider the notations $\Phi_{\snula}$ as very inconvenient. However, since it is commonly used, I decided to follow it.
\par
The family of all open $\varphi$-covers of a~topological space $X$ is usually denoted as $\Phi$ eventually as~$\Phi(X)$. However, the family of all Borel \mbox{$\varphi$-covers}  of a~topological space $X$ is usually denoted as ${\mathcal B}_{\Phi}$ eventually ${\mathcal B}_{\Phi(X)}$. I have decided for the notation $\Phi(\Borel)$. I~suppose that is more natural and can be used in more general cases as we already did.
%%%%%%%%%%%%%%%%%%%%%%%%%%%%%%%%%%%%%%%%%%%%%%%%%%
%%%%%%
%%%%%%%%%%%%%%%%%%%%%%%%%%%%%%%%%%%%%%%%%%%%%%%%%%%%%%%
%%%%%%%%%%%%%%%%%%%%%%%%%%%%%%%%%%%%%%%%%%%%%%%%%%%%%%%

%%%%%%%%%%%%%%%%%%%%%%%%%%%%%%%%%%%%%%%%%%%%%%%%%%%%%%%
\end{document}